\documentclass{article}

\usepackage[T1]{fontenc}
\usepackage[utf8]{inputenc}
\usepackage{setspace}
\usepackage{authblk}
\usepackage{float}
\usepackage{cuted}
\usepackage{booktabs}
\usepackage{multirow}
\usepackage{amsmath,amssymb}
\usepackage{graphicx}
\usepackage{array}
\usepackage{comment}
\usepackage{xcolor}
\usepackage{psfrag}
\usepackage{enumerate}
\usepackage{caption}
\usepackage{graphicx}
\usepackage{url}
\usepackage[natbibapa]{apacite}

\usepackage{url}
\usepackage[colorlinks=true,linkcolor=blue,citecolor=blue,urlcolor=blue]{hyperref}
\usepackage{doi}

\title{\bf An Adaptive Upper One-Sided Cumulative Sum Control Chart with Joint Parameter Optimization for Monitoring the Ratio of Two Normal Variables in Short Production Runs}

\author[1,2]{Kim Duc Tran}

\affil[1]{International Chair in DS \& XAI, International Research Institute for Artificial Intelligence and Data Science, Dong A University, Danang, Vietnam}

\affil[2]{Univ. Lille, ENSAIT, ULR 2461 - GEMTEX - G\'enie et Mat\'eriaux Textiles, F-59000 Lille, France}

\begin{document}
\onehalfspacing

\maketitle
\footnotetext[0]{Corresponding author: Kim Duc Tran, Email: ductk@donga.edu.vn}

\begin{abstract}
Monitoring the ratio of two correlated normal variables is increasingly important in statistical process control, since many quality characteristics are expressed in relative rather than absolute form. Memory-type ratio charts have mostly been developed for long production runs, while their finite-horizon counterparts rely on a fixed reference value $k$ derived from a specified shift. Such fixed-$k$ designs are not optimal at a given out-of-control magnitude and, in low-variability regimes, yield \emph{boundary solutions} for which the in-control truncated average run length (TARL$_0$) is unattainable.

This paper proposes an upper one-sided cumulative sum (CUSUM) control chart for the ratio $Z=X/Y$ in short production runs, denoted CUSUM-RZ$^+$ (RZ standing for the ratio $Z$), with fully adaptive joint optimization of $k$ and the decision interval $h$. Given a target TARL$_0=I$ and shift $\tau_{\mathrm{target}}$, a bilevel problem calibrates $h(k)$ by inner root-finding to satisfy the TARL$_0$ constraint and selects $k^\star$ by outer line search to minimize the out-of-control TARL$_1$. Both use a finite-state Markov-chain framework with an accurate ratio approximation; the inner step recovers boundary cases that fixed-$k$ designs cannot.

The chart is assessed through matched-horizon benchmarks against Shewhart-RZ, exponentially weighted moving average (EWMA-RZ), and fixed-$k$ CUSUM-RZ$^+$ charts, Monte Carlo robustness studies, and a Phase~I estimation analysis. All memory-type charts outperform the Shewhart-RZ baseline; the adaptive design matches them under stable correlation and improves appreciably when correlation rises from Phase~I to Phase~II. It is insensitive to symmetric heavy tails yet mildly anti-conservative under contamination, and $m\geq 100$ subgroups keep the TARL$_0$ relative bias near $1\%$.
\end{abstract}


\section{Introduction}
\label{sec:introduction}

Statistical process control (SPC) provides a fundamental set of tools for monitoring process stability and detecting abnormal changes in key quality characteristics. Among the most widely used SPC procedures are the Shewhart, cumulative sum (CUSUM), and exponentially weighted moving average (EWMA) control charts, with memory-type charts being well-known for their superior sensitivity to small and moderate process shifts compared with classical Shewhart designs \citep{Montgomery2019}.

In many practical applications, the quality characteristic of interest is not a single process variable but the ratio between two related variables. Examples include strength-to-weight indices, yield or content ratios, concentration measures, and several efficiency indicators encountered in manufacturing, process industries, and service systems. Monitoring such ratio-type characteristics is statistically challenging because the ratio of two normal variables does not follow a standard normal distribution and is typically asymmetric, with a shape that depends on the means, variances, and correlation of the underlying variables \citep{Cedilnik2004,Faraz2014,PhamGia2006}. For this reason, dedicated ratio-specific control charts are needed.

The literature on ratio monitoring has expanded substantially over the last decade. Early studies focused on deriving accurate approximations to the distribution of the ratio statistic and on developing Shewhart-type ratio charts for Phase~II monitoring \citep{CelanoCastagliola2016a,CelanoCastagliola2016b,Faraz2014}. These ideas were later extended to more advanced memory-type schemes, including run-rules, EWMA, and CUSUM-type charts, which were shown to provide improved detection capability when the process ratio experiences persistent but moderate shifts \citep{TranCastagliolaCelano2016,TranRunRules2016,CelanoCastagliolaTran2018CUSUM}. Most of these contributions, however, were developed under the assumption of long production runs, where the monitoring horizon is effectively infinite and chart performance is assessed through the traditional average run length (ARL).

This long-run assumption is increasingly unrealistic in modern manufacturing environments. High-mix low-volume production, customized manufacturing, job-shop systems, pilot production, and start-up phases often generate only a limited number of samples during a production campaign. In such settings, the process terminates after a finite number of inspections, and the classical ARL becomes operationally less meaningful. Consequently, truncated performance measures such as the truncated average run length (TARL) are more appropriate for evaluating control-chart performance in short production runs \citep{Hu2024,Tran2021,Yang2025}.

Recent studies have begun to address ratio monitoring in finite-horizon settings. \citet{Tran2021} proposed one-sided Shewhart charts for monitoring the ratio of two normal variables in short production runs and used TARL as the primary design criterion. \citet{Hu2024} extended this work to memory-based schemes, including EWMA- and CUSUM-type charts, using a Markov-chain framework. \citet{Yang2025} further investigated fixed and variable sampling-interval strategies for CUSUM-type ratio monitoring in short production runs. These contributions clearly establish the feasibility of short-run monitoring for ratio-type characteristics.

Nevertheless, an important methodological limitation of the existing finite-horizon CUSUM-RZ literature has not been fully addressed. In all of these works, the reference value $k$ is fixed in advance from a target shift size (often $\tau_{\mathrm{ref}}=1.05$, leading to $k=1.025$ when $z_0=1$), while only the decision interval $h$ is optimized through the TARL$_0$ constraint. This convention is inherited from the long-run CUSUM literature, where the reference value is asymptotically optimal at the targeted shift; however, in short production runs it has two consequences that have been observed but not fully resolved. First, the fixed-$k$ design is not generally optimal at any specific out-of-control magnitude $\tau_{\mathrm{target}}$, because the optimal reference value depends on the joint distribution of the subgroup ratio statistic, which itself depends on the coefficients of variation, the correlation structure, and the subgroup size. Second, when the in-control coefficients of variation are extremely small (e.g., $\gamma_X=\gamma_Y=0.01$), the fixed reference $k=1.025$ is far above the in-control mean of the ratio statistic, so that the upward signal probability is too small for the prescribed TARL$_0$ to be attainable within the finite horizon. Numerically, the calibration then collapses to a lower bound of $h$, producing what we term \emph{boundary solutions}.

Motivated by these limitations, this paper proposes an upper one-sided cumulative sum chart for monitoring increases in the ratio of two correlated normal variables in short production runs, with a fully adaptive joint optimization of $(k,h)$. The proposed chart, denoted by CUSUM-RZ$^+$ (adaptive), is calibrated by solving a constrained bilevel design problem: for any candidate reference value $k$, an inner one-dimensional root-finding step computes the corresponding $h(k)$ so that the in-control TARL$_0$ matches the prescribed target $I$, and an outer line search over a grid of admissible $k$ values selects the $k^\star$ that minimizes the out-of-control TARL$_1$ at a user-specified target shift $\tau_{\mathrm{target}}$. Both steps rely on the Markov-chain framework of \citet{BrookEvans1972} adapted to a finite truncation horizon.

The contributions of the paper are fourfold. First, the adaptive joint $(k,h)$ design provides a transparent and computationally efficient alternative to fixed-$k$ short-run CUSUM charts, with a clear interpretation in the bivariate normal ratio setting. Second, the inner $h$-calibration step is shown to systematically recover the boundary cases of fixed-$k$ designs by simply moving $k$ closer to the in-control ratio $z_0$, thereby ensuring numerical feasibility across a wide range of variability regimes. Third, the chart is benchmarked under matched finite-horizon conditions against three competing short-run ratio charts (Shewhart-RZ, EWMA-RZ, and the fixed-$k$ CUSUM-RZ$^+$), establishing its competitive performance in regimes of practical interest, with particularly large gains over the Shewhart-RZ baseline. Fourth, the paper documents the robustness of the adaptive chart under departures from joint normality and the impact of Phase~I parameter estimation on the actually achieved TARL$_0$ — two practical aspects that are rarely investigated in the short-run ratio-monitoring literature.

The remainder of the paper is organized as follows. Section~\ref{sec:litreview} reviews the relevant literature and clarifies the research position relative to recent short-run ratio charts. Section~\ref{sec:distribution} presents the statistical model and the distributional approximation of the ratio statistic, together with a Monte Carlo validation. Section~\ref{sec:implementation} introduces the proposed upper one-sided CUSUM-RZ$^+$ chart and its statistical properties. Section~\ref{sec:adaptive_design} describes the adaptive joint $(k,h)$ design via the constrained bilevel approach. Section~\ref{sec:numerical} reports the numerical analysis, including design tables, comparisons with the fixed-$k$ design, matched-horizon benchmarks, robustness diagnostics, and Phase~I results. Section~\ref{sec:illustrative} illustrates the chart through a food-industry example. The last section concludes and outlines directions for future research.

\section{Literature review}
\label{sec:litreview}

The literature relevant to the present study can be organized into three closely related streams: (i) the statistical theory of the ratio of two normal variables; (ii) control charts developed specifically for ratio-type quality characteristics; and (iii) monitoring schemes designed for short production runs or finite-horizon processes. Together, these streams provide the theoretical and methodological foundation for the development of the adaptive CUSUM-RZ$^+$ chart.

\subsection{Ratio of two normal variables}

The ratio of two normal variables has been studied for many decades because of its mathematical difficulty and practical importance. A seminal contribution is due to \citet{Geary1930}, who investigated the quotient of two normal variates and highlighted the nonstandard nature of the resulting distribution. \citet{Fieller1932} and later \citet{Hinkley1969} provided exact distributional results under joint normality. Subsequent studies refined this line of work by providing approximations and analytical characterizations under different dependence structures, including \citet{Hayya1975}, \citet{Cedilnik2004}, and \citet{PhamGia2006}.

These studies established that the ratio distribution is generally asymmetric, may exhibit heavy tails, and depends on the means, variances, and correlation coefficient of the underlying variables. Such properties make ratio monitoring fundamentally different from monitoring a single normal quality characteristic. In many industrial contexts, moreover, the ratio itself is the parameter of engineering interest, for example in strength-to-weight, output-to-input, fill-to-capacity, and concentration-type indices \citep{Faraz2014,Oksoy1994,Spisak1990}. These applications motivated the transition from purely distributional studies to dedicated SPC procedures for directly monitoring the ratio statistic.

\subsection{Control charts for monitoring the ratio of two normal variables}

A second stream of research concerns control charts designed specifically for the ratio of two normal variables, denoted by $RZ$. Early work includes the ratio chart of \citet{Spisak1990} and the quotient-based SPC formulation of \citet{Oksoy1994}. More systematic developments appeared with the work of \citet{Faraz2014}, who investigated the statistical performance of a Shewhart-type chart for individual observations monitoring the ratio of two normal variables. Subsequent studies by \citet{CelanoCastagliola2016a,CelanoCastagliola2016b} developed Phase~II Shewhart and synthetic charts for the same problem, showing that one-sided designs are often preferable because of the asymmetry of the ratio distribution.

To improve sensitivity to small and moderate shifts, \citet{TranRunRules2016} proposed run-rules-type charts, and \citet{TranCastagliolaCelano2016} introduced EWMA-type charts for monitoring the ratio of two normal variables. \citet{CelanoCastagliolaTran2018CUSUM} developed one-sided CUSUM charts for monitoring the ratio of population means of a bivariate normal distribution, demonstrating their strong sensitivity to persistent shifts in the ratio parameter. Methodological refinements based on auxiliary information were later reported by \citet{NguyenAux2019}. These studies established the long-run Phase~II foundation for CUSUM-type ratio monitoring on which the present paper builds.

\subsection{Short-run and finite-horizon monitoring schemes}

The third stream concerns short-run statistical process control. Classical SPC was originally developed for mass-production settings, but many modern manufacturing environments operate under high product variety, frequent setup changes, and limited production quantities. In such settings, only a small number of samples may be available during a given production run, and long-run metrics such as the ARL become less informative \citep{DelCastillo1996,Farnum1992,Lin2005}.

A more formal treatment of short-run monitoring emerged through finite-horizon and truncated-run-length formulations. \citet{CelanoEtAl2011t} developed Shewhart and EWMA $t$ control charts for short production runs, which remain among the most widely studied finite-horizon schemes for monitoring a single normal characteristic with unknown variance. More recent studies extended truncated-performance ideas to other parameters, including the multivariate coefficient of variation \citep{MCVSPR2023,FiniteHorizonScale2024}.

The relevance of this short-run literature to ratio monitoring is direct. \citet{Tran2021} proposed one-sided Shewhart charts for monitoring the ratio of two normal variables in short production runs and used TARL as the main design criterion. \citet{Hu2024} later developed memory-based EWMA- and CUSUM-type schemes for the same problem using a Markov-chain approximation. \citet{Yang2025} extended finite-horizon ratio monitoring by incorporating fixed and variable sampling-interval strategies into CUSUM-type designs. Taken together, these works establish the feasibility of finite-horizon ratio monitoring, yet they all fix the reference value $k$ in advance and calibrate only the decision interval $h$; the joint design of $(k,h)$ — and the boundary infeasibility it induces when the coefficients of variation are very small — remains unaddressed.

\subsection{Research position: adaptive joint $(k,h)$ design and boundary recovery}
\label{sec:research_position}

In all of the finite-horizon CUSUM-RZ designs cited above, the reference value $k$ is fixed by a target shift through $k=z_0(1+\tau_{\mathrm{ref}})/2$ and only $h$ is calibrated to satisfy the TARL$_0$ constraint. While this design is computationally simple and inherits the asymptotic optimality of the long-run CUSUM at the chosen $\tau_{\mathrm{ref}}$, two practical limitations have been observed.

First, the fixed-$k$ design is not jointly optimal in finite-horizon settings, because the optimal reference value depends on the joint distribution of $\hat Z_i=\bar X_i/\bar Y_i$ — and therefore on $(\gamma_X,\gamma_Y,\rho_0,n,I)$ — rather than only on the targeted shift. As a result, for a given out-of-control magnitude $\tau_{\mathrm{target}}$, the CUSUM-RZ$^+$ chart with fixed $k=1.025$ is dominated, often substantially, by a chart whose $k$ is selected jointly with $h$ under the TARL$_0$ constraint. Second, when the coefficients of variation are very small, the in-control distribution of $\hat Z_i$ is concentrated tightly around $z_0=1$, while the fixed reference $k=1.025$ remains far above this mean. The increment $\hat Z_i-k$ is then predominantly negative under the in-control state and the upward signal probability becomes vanishingly small, so that the calibration equation $\mathrm{TARL}_0(h)=\tau_0$ cannot be satisfied within the admissible range $h\geq h_{\min}$. Numerical procedures collapse to the lower bound, producing \emph{boundary solutions} that are not exact TARL$_0$ calibrations.

These observations motivate the present paper, which positions itself as a focused but distinct extension of the recent short-run ratio-monitoring literature in three directions:

\begin{enumerate}
\item \textbf{Adaptive joint $(k,h)$ design.} The reference value $k$ is no longer fixed; instead, a constrained bilevel procedure jointly determines $(k^\star,h^\star)$ that minimizes TARL$_1$ at a user-specified $\tau_{\mathrm{target}}$ subject to the TARL$_0$ constraint. This is closer in spirit to optimal CUSUM design \citep{MoustakidesPolunchenko2011,Reynolds2014} but is adapted to finite-horizon TARL.
\item \textbf{Boundary recovery.} Because the inner step calibrates $h(k)$ for any candidate $k$ in an admissible range that includes values close to $z_0$, the adaptive chart automatically resolves the boundary configurations of the fixed-$k$ design. This effectively eliminates the unattainable cases of the recent fixed-$k$ literature.
\item \textbf{Three-way diagnostic study.} The chart is benchmarked under matched finite-horizon conditions against three competing charts (Shewhart-RZ, EWMA-RZ, fixed-$k$ CUSUM-RZ$^+$); its robustness is assessed under non-normal models (Student-$t$ and contaminated bivariate normal); and the impact of Phase~I estimation on the achieved TARL$_0$ is quantified as a function of the calibration sample size.
\end{enumerate}

This combination of contributions has not, to the best of our knowledge, been jointly addressed in the short-run ratio-monitoring literature.

\section{The distribution of the ratio $Z$}
\label{sec:distribution}

Let $X$ and $Y$ be two correlated normal random variables such that
\[
\mathbf{W}=(X,Y)^T \sim N(\pmb{\mu}_{\mathbf{W}},\pmb{\Sigma}_{\mathbf{W}}),
\]
where $\mathbf{W}$ is a bivariate normal random vector with mean vector and variance--covariance matrix
\begin{align}
\pmb{\mu}_{\mathbf{W}} &=
\begin{pmatrix}
\mu_X\\
\mu_Y
\end{pmatrix},
&
\pmb{\Sigma}_{\mathbf{W}} &=
\begin{pmatrix}
\sigma_X^2 & \rho \sigma_X \sigma_Y\\
\rho \sigma_X \sigma_Y & \sigma_Y^2
\end{pmatrix}.
\label{eq:mu_sigma_W}
\end{align}
The coefficients of variation of $X$ and $Y$ are
\[
\gamma_X=\frac{\sigma_X}{\mu_X},
\qquad
\gamma_Y=\frac{\sigma_Y}{\mu_Y},
\]
and the standard-deviation ratio is $\omega=\sigma_X/\sigma_Y$.

Let $Z=X/Y$ denote the ratio of the two quality characteristics. Because the exact distribution of $Z$ is analytically cumbersome, an accurate approximation is commonly employed in the ratio-monitoring literature \citep{CelanoCastagliola2016a,TranCastagliolaCelano2016,CelanoCastagliolaTran2018CUSUM}. Specifically, the cumulative distribution function (c.d.f.) of $Z$ can be approximated by
\begin{align}
F_Z(z \mid \gamma_X,\gamma_Y,\omega,\rho)
\simeq
\Phi\!\left(\frac{A}{B}\right),
\label{eq:Zcdf}
\end{align}
where
\[
A=\frac{z}{\gamma_Y}-\frac{\omega}{\gamma_X},
\qquad
B=\sqrt{\omega^2-2\rho\omega z+z^2},
\]
and $\Phi(\cdot)$ denotes the c.d.f.\ of the standard normal distribution. Differentiating \eqref{eq:Zcdf} yields the approximate probability density function
\begin{align}
f_Z(z \mid \gamma_X,\gamma_Y,\omega,\rho)
\simeq
\left(
\frac{1}{B\gamma_Y}
-
\frac{(z-\rho\omega)A}{B^3}
\right)
\phi\!\left(\frac{A}{B}\right),
\label{eq:Zpdf}
\end{align}
where $\phi(\cdot)$ is the standard normal p.d.f.

An approximation of the inverse distribution function (i.d.f.) is obtained by solving $F_Z(z\mid\cdot)=p$ for $0<p<1$, which yields the closed-form expression
\begin{align}
F_Z^{-1}(p \mid \gamma_X,\gamma_Y,\omega,\rho)
\simeq
\begin{cases}
\dfrac{-C_2-\sqrt{C_2^2-4C_1C_3}}{2C_1}, & \text{if } p\in(0,0.5],\\[2ex]
\dfrac{-C_2+\sqrt{C_2^2-4C_1C_3}}{2C_1}, & \text{if } p\in[0.5,1),
\end{cases}
\label{eq:Zidf}
\end{align}
with
\begin{align*}
C_1 &= \frac{1}{\gamma_Y^2}-\left(\Phi^{-1}(p)\right)^2,\\
C_2 &= 2\omega\left[\rho\left(\Phi^{-1}(p)\right)^2-\frac{1}{\gamma_X\gamma_Y}\right],\\
C_3 &= \omega^2\left[\frac{1}{\gamma_X^2}-\left(\Phi^{-1}(p)\right)^2\right].
\end{align*}
The approximations in \eqref{eq:Zcdf}--\eqref{eq:Zidf} are central to the construction of ratio-type control charts because they allow probability limits, transition probabilities, and run-length characteristics to be evaluated efficiently.

\subsection{Validation of the distributional approximation}
\label{sec:approx_validation}

Because the transition probabilities used in the Markov-chain design are computed from the c.d.f.\ and i.d.f.\ of the subgroup ratio statistic $\hat Z_i=\bar X_i/\bar Y_i$, it is important to assess the numerical adequacy of the approximation. To that end, we conducted a Monte Carlo validation study based on $150{,}000$ replications for each representative parameter configuration. Bivariate normal subgroup means were generated under the corresponding values of $n$, $(\gamma_X,\gamma_Y)$, and $\rho$, and the empirical distribution of $\hat Z_i$ was compared with the analytical approximation through (i) the Kolmogorov--Smirnov distance between the empirical and approximated c.d.f.s, and (ii) the signed errors in the $0.05$, $0.50$, and $0.95$ quantiles, $\Delta q_p=q_p^{(\mathrm{Approx})}-q_p^{(\mathrm{MC})}$.

The results are reported in Table~\ref{tab:approx_validation}. Across all settings, the Kolmogorov--Smirnov distance remains below $0.003$, and the absolute quantile error never exceeds $0.0013$ at the three reported probability levels. These results indicate that the approximation is sufficiently accurate for the present calibration and performance-evaluation purposes, including settings with small and large coefficients of variation, unequal coefficients of variation, and both negative and positive correlation structures.

\begin{table}[!htbp]
\centering
\caption{Monte Carlo validation of the approximation used for the distribution of the subgroup ratio statistic $\hat Z_i$. The column $\Delta q_p$ reports the signed quantile error $q_p^{(\mathrm{Approx})}-q_p^{(\mathrm{MC})}$.}
\label{tab:approx_validation}
\setlength{\tabcolsep}{4pt}
\renewcommand{\arraystretch}{1.05}
\begin{tabular}{ccccccc}
\toprule
$n$ & $(\gamma_X,\gamma_Y)$ & $\rho$ & KS distance & $\Delta q_{0.05}$ & $\Delta q_{0.50}$ & $\Delta q_{0.95}$ \\
\midrule
1  & $(0.01,0.01)$ & $-0.8$ & 0.0029 &  0.0001 & -0.0001 &  0.0001 \\
1  & $(0.20,0.20)$ &  $0.8$ & 0.0029 & -0.0012 & -0.0005 & -0.0012 \\
5  & $(0.20,0.20)$ &  $0.8$ & 0.0016 &  0.0005 &  0.0000 &  0.0002 \\
5  & $(0.01,0.20)$ &  $0.4$ & 0.0017 &  0.0003 &  0.0000 &  0.0003 \\
15 & $(0.20,0.01)$ & $-0.4$ & 0.0013 & -0.0005 & -0.0001 &  0.0001 \\
\bottomrule
\end{tabular}
\end{table}

\section{Implementation of the upper one-sided CUSUM-RZ$^+$ chart}
\label{sec:implementation}

Consider a short production run in which a lot of $N$ items is produced over a fixed rolling horizon $H$. Let $I$ denote the number of planned inspections during this horizon, so that the time interval between two consecutive inspections is $\mathcal{S}_h=H/(I+1)$. At each inspection epoch, a sample of size $n$ is collected and the bivariate quality characteristic $\mathbf{W}$ is measured for each sampled unit. Let
\[
\mathbf{W}_{i,j}=(X_{i,j},Y_{i,j})^T,
\qquad
j=1,2,\ldots,n,
\]
denote the $j$-th observation collected at the $i$-th inspection. We assume that
\[
\mathbf{W}_{i,j}\sim N(\pmb{\mu}_{\mathbf{W},i},\pmb{\Sigma}_{\mathbf{W},i}),
\qquad i=1,2,\ldots,I,
\]
with
\begin{align}
\pmb{\mu}_{\mathbf{W},i} &=
\begin{pmatrix}
\mu_{X,i}\\
\mu_{Y,i}
\end{pmatrix},
&
\pmb{\Sigma}_{\mathbf{W},i} &=
\begin{pmatrix}
\sigma_{X,i}^2 & \rho\,\sigma_{X,i}\sigma_{Y,i}\\
\rho\,\sigma_{X,i}\sigma_{Y,i} & \sigma_{Y,i}^2
\end{pmatrix}.
\label{eq:Wi_distribution}
\end{align}

To design the proposed chart, let $\gamma_X$ and $\gamma_Y$ be known and constant coefficients of variation, and let $z_0=\mu_{X,i}/\mu_{Y,i}$ and $\rho_0$ denote the in-control values of the ratio and correlation. Following common practice in the short-run ratio-monitoring literature, we assume that the standard deviations are proportional to the corresponding means, i.e.\ $\sigma_{X,i}=\gamma_X \mu_{X,i}$ and $\sigma_{Y,i}=\gamma_Y \mu_{Y,i}$, so that the coefficients of variation are preserved across $i$.

The monitoring statistic at the $i$-th inspection is the ratio of sample means,
\begin{align}
\hat Z_i
=
\frac{\bar X_i}{\bar Y_i}
=
\frac{\sum_{j=1}^{n}X_{i,j}}{\sum_{j=1}^{n}Y_{i,j}},
\qquad i=1,2,\ldots,I.
\label{eq:ratio_statistic}
\end{align}
Since $\bar X_i\sim N(\mu_{X,i},\sigma_{X,i}^2/n)$ and $\bar Y_i\sim N(\mu_{Y,i},\sigma_{Y,i}^2/n)$, the coefficients of variation of the sample means are $\gamma_{\bar X}=\gamma_X/\sqrt{n}$ and $\gamma_{\bar Y}=\gamma_Y/\sqrt{n}$, while the in-control standard-deviation ratio is
\begin{align}
\omega_0=z_0\frac{\gamma_X}{\gamma_Y}.
\label{eq:omega_0}
\end{align}

Using the approximation developed in Section~\ref{sec:distribution}, the c.d.f.\ and i.d.f.\ of $\hat Z_i$ under the in-control state are
\begin{align}
F_{\hat Z_i}(z \mid n,\gamma_X,\gamma_Y,z_0,\rho_0)
&=
F_Z\!\left(
z \,\middle|\, \frac{\gamma_X}{\sqrt{n}},
\frac{\gamma_Y}{\sqrt{n}},
\frac{z_0\gamma_X}{\gamma_Y},
\rho_0
\right),
\label{eq:FhatZ_ic}
\\
F_{\hat Z_i}^{-1}(p \mid n,\gamma_X,\gamma_Y,z_0,\rho_0)
&=
F_Z^{-1}\!\left(
p \,\middle|\, \frac{\gamma_X}{\sqrt{n}},
\frac{\gamma_Y}{\sqrt{n}},
\frac{z_0\gamma_X}{\gamma_Y},
\rho_0
\right).
\label{eq:FhatZinv_ic}
\end{align}
Under an out-of-control state characterized by a shifted ratio $z_1=\tau z_0$ (and possibly a different correlation $\rho_1$), the same approximation is used after replacing $(z_0,\rho_0)$ with $(z_1,\rho_1)$.

The proposed upper one-sided CUSUM-RZ$^+$ statistic is defined as
\begin{align}
C_0^+ = 0,
\qquad
C_i^+ = \max\left\{0,\, C_{i-1}^+ + (\hat Z_i-k)\right\},
\qquad i=1,2,\ldots,I,
\label{eq:cusum_upper}
\end{align}
where $k$ is the reference value. An out-of-control signal is triggered as soon as
\begin{align}
C_i^+ \ge h,
\label{eq:signal_upper}
\end{align}
where $h>0$ is the decision interval. Since the present paper focuses on detecting upward shifts only, no lower-sided chart is considered.

In the long-run CUSUM literature, $k$ is typically chosen as the midpoint between $z_0$ and a target out-of-control value $z_1=\tau_{\mathrm{ref}}z_0$, that is,
\begin{align}
k_{\mathrm{fixed}}=\frac{z_0+\tau_{\mathrm{ref}}z_0}{2}=z_0\frac{1+\tau_{\mathrm{ref}}}{2},
\label{eq:k_reference_fixed}
\end{align}
with $\tau_{\mathrm{ref}}=1.05$ giving $k_{\mathrm{fixed}}=1.025$ when $z_0=1$. As discussed in Section~\ref{sec:research_position}, this fixed-$k$ choice is suboptimal in finite-horizon settings and produces boundary solutions when the coefficients of variation are very small. Section~\ref{sec:adaptive_design} therefore introduces an adaptive joint $(k,h)$ design that addresses both issues.

\section{Adaptive joint $(k,h)$ design via constrained bilevel optimization}
\label{sec:adaptive_design}

This section presents the proposed adaptive design. We first review the truncated run length and the Markov-chain framework, then introduce the bilevel optimization algorithm that yields $(k^\star,h^\star)$.

\subsection{Truncated run length and TARL}

Let $RL$ denote the run length of the upper one-sided CUSUM-RZ$^+$ chart. Since the production run contains only $I$ planned inspections, the run length is naturally truncated at the end of the horizon, giving
\begin{align}
RL^{(T)}=\min(RL,I+1),
\label{eq:trl}
\end{align}
where $I+1$ corresponds to the event that no signal is produced before the end of the run. The in-control and out-of-control TARLs are
\begin{align}
\mathrm{TARL}_0 &= \mathbb{E}_0\!\left[RL^{(T)}\right],
&
\mathrm{TARL}_1 &= \mathbb{E}_1\!\left[RL^{(T)}\right],
\label{eq:tarl_def}
\end{align}
where the subscripts refer to the in-control and out-of-control states, respectively.

\subsection{Markov-chain framework}
\label{sec:markov_chain}

Following \citet{BrookEvans1972} and the recent short-run ratio-monitoring literature \citep{Hu2024,Yang2025}, the continuation region $[0,h)$ is partitioned into $m$ equal subintervals of width $\Delta=h/m$, with representative midpoints
\begin{align}
c_r=\left(r-\frac{1}{2}\right)\Delta,
\qquad r=1,2,\ldots,m.
\end{align}
Suppose the current state is $c_r$. Then by \eqref{eq:cusum_upper} the next CUSUM statistic is $C_i^+=\max\{0,c_r+\hat Z_i-k\}$. Let $F_{\hat Z}^{(j)}$ denote the c.d.f.\ of $\hat Z_i$ under state $j\in\{0,1\}$. The transition probability from state $r$ to state $s$ is
\begin{align}
p_{rs}^{(j)}=
\begin{cases}
F_{\hat Z}^{(j)}\!\left(k-c_r+\Delta\right), & s=1,\\[1.25ex]
F_{\hat Z}^{(j)}\!\left(k-c_r+s\Delta\right)
-
F_{\hat Z}^{(j)}\!\left(k-c_r+(s-1)\Delta\right), & s=2,\ldots,m,
\end{cases}
\label{eq:trans_upper}
\end{align}
and the absorbing probability is $p_{rA}^{(j)}=1-F_{\hat Z}^{(j)}(k-c_r+h)$. Let $\mathbf{Q}^{(j)}=[p_{rs}^{(j)}]$ be the $m\times m$ transient transition matrix.

Let $\boldsymbol{\xi}$ denote the initial probability vector of the Markov chain. Since $C_0^+=0$, the initial probability mass is assigned to the first subinterval. Then for any $j\in\{0,1\}$,
\begin{align}
\mathrm{TARL}_j(k,h)
&=
\sum_{i=0}^{I}
\boldsymbol{\xi}\left(\mathbf{Q}^{(j)}\right)^i\mathbf{1}
\nonumber\\
&=
\boldsymbol{\xi}
\left[
\mathbf{I}-\left(\mathbf{Q}^{(j)}\right)^{I+1}
\right]
\left[
\mathbf{I}-\mathbf{Q}^{(j)}
\right]^{-1}
\mathbf{1},
\label{eq:tarl_matrix}
\end{align}
where $\mathbf{1}$ is a column vector of ones. The matrix form makes explicit the dependence of $\mathrm{TARL}_j$ on the design pair $(k,h)$ through both $\mathbf{Q}^{(j)}$ and the absorbing structure. Throughout the paper we use $m=60$ states, which is sufficient for stable evaluation in the parameter ranges considered here.

\subsection{Constrained bilevel design problem}
\label{sec:bilevel}

Given an in-control TARL$_0$ target $\tau_0=I$ and an out-of-control target shift $\tau_{\mathrm{target}}>1$, the adaptive design problem is
\begin{align}
(k^\star,h^\star)
&=
\arg\min_{k,h}\;
\mathrm{TARL}_1(k,h;\tau_{\mathrm{target}})
\nonumber\\
&\text{subject to}\quad
\mathrm{TARL}_0(k,h)=\tau_0,
\qquad
k\in\mathcal{K},\;
h\in[h_{\min},h_{\max}],
\label{eq:bilevel}
\end{align}
where $\mathcal{K}=[z_0,\,z_0\tau_{\max}]$ is the admissible range of reference values, with $\tau_{\max}>\tau_{\mathrm{target}}$ to allow $k^\star$ to lie strictly between $z_0$ and the upper bound. In the present paper we use $\mathcal{K}=[1.000,1.100]$ and $[h_{\min},h_{\max}]=[10^{-3},10]$.

The problem in \eqref{eq:bilevel} is a constrained bilevel optimization. We solve it by exploiting two facts: (i) for any fixed $k\in\mathcal{K}$, $\mathrm{TARL}_0(k,h)$ is monotonically increasing in $h$ within $[h_{\min},h_{\max}]$; and (ii) the outer objective $\mathrm{TARL}_1(k,h(k);\tau_{\mathrm{target}})$ is empirically unimodal in $k$ over $\mathcal{K}$ for the parameter ranges considered here. The first property guarantees that for any $k$, a unique $h(k)$ satisfies the TARL$_0$ constraint, which can be located by bisection or Brent's method. The second property motivates a one-dimensional grid search over $\mathcal{K}$ followed by a local refinement around the best grid point.
Algorithm~\ref{alg:bilevel} summarizes the procedure.

\begin{table}[!htbp]
\centering
\caption{Adaptive joint $(k,h)$ design via constrained bilevel optimization (pseudo-code).}
\label{alg:bilevel}
\setlength{\tabcolsep}{4pt}
\renewcommand{\arraystretch}{1.10}
\resizebox{\columnwidth}{!}{%
\begin{tabular}{l}
\toprule
\textbf{Algorithm 1: Adaptive design for the upper one-sided CUSUM-RZ$^+$ chart} \\
\midrule
\textbf{Input:} $z_0$, $\tau_{\mathrm{target}}$, $I$, $n$, $\gamma_X$, $\gamma_Y$, $\rho_0$, $\rho_1$, $m=60$, grid $\mathcal{K}_g$ over $\mathcal{K}$\\
\textbf{Output:} $(k^\star,h^\star)$, $\mathrm{TARL}_0(k^\star,h^\star)$, $\mathrm{TARL}_1(k^\star,h^\star;\tau_{\mathrm{target}})$\\
\\
\textbf{1. Inner step (constraint satisfaction).} For each candidate $k\in\mathcal{K}_g$:\\
\quad (a) Define $g_k(h)=\mathrm{TARL}_0(k,h)-\tau_0$ on $[h_{\min},h_{\max}]$.\\
\quad (b) Locate $h(k)$ by Brent's method; if $g_k(h_{\min})\geq 0$, set $h(k)=h_{\min}$ (boundary).\\
\\
\textbf{2. Outer step (objective minimization).} For each $k\in\mathcal{K}_g$:\\
\quad (a) Compute $\mathbf{Q}^{(1)}$ with shifted parameters $(z_1,\rho_1)$, $z_1=\tau_{\mathrm{target}}z_0$.\\
\quad (b) Compute $T_1(k)=\mathrm{TARL}_1(k,h(k);\tau_{\mathrm{target}})$ from \eqref{eq:tarl_matrix}.\\
\\
\textbf{3. Selection.} Set $k^\star=\arg\min_{k\in\mathcal{K}_g} T_1(k)$ and $h^\star=h(k^\star)$.\\
\\
\textbf{4. Diagnostics.} Return feasibility flag (whether $\mathrm{TARL}_0(k^\star,h^\star)$ matches $\tau_0$\\
\quad within numerical tolerance) and boundary flag (whether $h^\star=h_{\min}$).\\
\bottomrule
\end{tabular}}
\end{table}

\subsection{Boundary recovery and feasibility}
\label{sec:boundary_recovery}

The bilevel structure of Algorithm~1 implies that the adaptive chart is feasible (i.e., the prescribed TARL$_0$ is exactly attained) whenever there exists $k\in\mathcal{K}$ such that $g_k(h_{\min})<0$. Because the inner search includes values of $k$ arbitrarily close to $z_0$, the upward signal probability under the in-control state can be made arbitrarily large by lowering $k$, and feasibility is guaranteed across the variability range considered in this paper. In the main numerical study (Section~\ref{sec:numerical}), which uses $(\gamma_X,\gamma_Y)\in\{(0.20,0.20),(0.01,0.20)\}$, all configurations admit a feasible interior solution and no boundary solution arises.

The boundary difficulty is confined to the extreme low-variability regime $\gamma_X=\gamma_Y=0.01$. There the fixed-$k$ design with $k=1.025$ produces an in-control increment $\hat Z_i-k$ that is almost surely negative, because $k=1.025$ lies far above the tightly concentrated in-control mean $z_0=1$. The statistic therefore never accumulates and the chart essentially never signals, so the achieved TARL$_0$ saturates at the truncation ceiling $I+1$ for \emph{every} $h\ge h_{\min}$, and this concentration only tightens as $n$ grows. The prescribed target $\tau_0<I+1$ is thus unattainable from above, the root-finder falls back to $h_{\min}$, and the boundary flag is raised. The adaptive design resolves this by lowering $k^\star$ toward $z_0$, which restores a nonzero upward signal probability and allows a finite $h^\star$ to attain $\tau_0$ exactly. Table~\ref{tab:R_boundary} reports this recovery for $n\in\{5,10,15\}$: in every case the fixed-$k$ design is infeasible (its decision interval is pinned at $h_{\min}=10^{-3}$ and the achieved TARL$_0$ equals the ceiling $I+1=31\neq 30$), whereas the adaptive design returns an interior $(k^\star,h^\star)$ that attains TARL$_0=30$.

\begin{table}[!htbp]
\centering
\caption{Boundary recovery in the extreme low-variability regime $\gamma_X=\gamma_Y=0.01$ ($z_0=1$, $\rho_0=0$, $\tau_{\mathrm{target}}=1.05$, $\tau_0=I=30$, $m=60$). The fixed-$k$ design cannot attain the prescribed TARL$_0$; the adaptive design lowers $k^\star$ toward $z_0$ and attains it exactly.}
\label{tab:R_boundary}
\setlength{\tabcolsep}{5pt}
\renewcommand{\arraystretch}{1.10}
\begin{tabular}{ccccccccc}
\toprule
\multirow{2}{*}{$n$} & \multicolumn{4}{c}{Fixed-$k$ design} & \multicolumn{4}{c}{Adaptive design}\\
\cmidrule(lr){2-5}\cmidrule(lr){6-9}
 & $k$ & $h$ & TARL$_0$ & feasible & $k^\star$ & $h^\star$ & TARL$_0$ & feasible\\
\midrule
$5$  & $1.0250$ & $0.0010^{\dagger}$ & $30.99$ & no  & $1.0161$ & $0.0021$ & $30.00$ & yes\\
$10$ & $1.0250$ & $0.0010^{\dagger}$ & $31.00$ & no  & $1.0083$ & $0.0047$ & $30.01$ & yes\\
$15$ & $1.0250$ & $0.0010^{\dagger}$ & $31.00$ & no  & $1.0064$ & $0.0043$ & $30.00$ & yes\\
\bottomrule
\end{tabular}

\vspace{3pt}
{\footnotesize $^{\dagger}$Pinned at the lower bound $h_{\min}=10^{-3}$; no $h\ge h_{\min}$ yields TARL$_0=30$ for the fixed-$k$ design.}
\end{table}

\subsection{Diagnostic interpretation of the optimal $k^\star$}
\label{sec:diagnostic_kstar}

\begin{figure}[H]
\centering
\includegraphics[width=0.95\columnwidth]{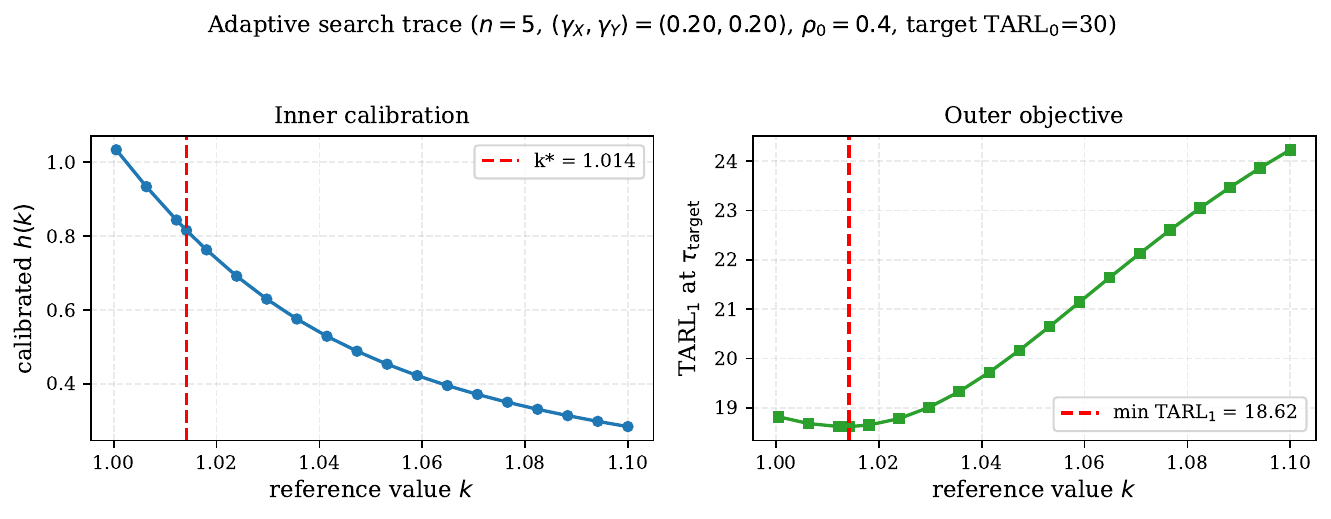}
\caption{Diagnostic of the bilevel adaptive search for $n=5$, $(\gamma_X,\gamma_Y)=(0.20,0.20)$, $\rho_0=0.4$, $\tau_{\mathrm{target}}=1.05$, and $\tau_0=I=30$. Left: inner-step calibration $h(k)$, monotonically decreasing in $k$. Right: outer-step objective $\mathrm{TARL}_1(k)$ at $\tau_{\mathrm{target}}$ as a function of $k$, with a clear U-shape and minimizer $k^\star\approx1.014<k_{\mathrm{fixed}}=1.025$.}
\label{fig:adaptive_trace}
\end{figure}

Figure~\ref{fig:adaptive_trace} illustrates Algorithm~1 for a representative case ($n=5$, $\gamma_X=\gamma_Y=0.20$, $\rho_0=0.4$, $\tau_{\mathrm{target}}=1.05$, $\tau_0=I=30$). The left panel shows the inner solution $h(k)$, which is monotonically decreasing in $k$: as $k$ moves away from $z_0$, the in-control increment $\hat Z_i-k$ becomes more negative on average, so a smaller $h$ suffices to attain the same in-control TARL$_0$. The right panel shows the outer objective $T_1(k)$, which exhibits a clear U-shaped pattern. The minimizer $k^\star\approx 1.014$ lies between $z_0=1$ and the long-run prescription $k_{\mathrm{fixed}}=1.025$, yielding $\mathrm{TARL}_1=18.62$ at $\tau_{\mathrm{target}}=1.05$. This confirms that, in finite-horizon settings, the optimal reference value differs from the long-run prescription, and the difference is a function of all the design parameters $(\gamma_X,\gamma_Y,\rho_0,n,I)$.

\section{Numerical analysis}
\label{sec:numerical}

This section reports a comprehensive numerical study of the proposed adaptive CUSUM-RZ$^+$ chart. Throughout the section, we set $z_0=1$, $\tau_{\mathrm{target}}=1.05$, $\tau_0=I=30$, and consider subgroup sizes $n\in\{5,10\}$. The variability and correlation grid focuses on parameter combinations of practical relevance, namely $(\gamma_X,\gamma_Y)\in\{(0.20,0.20),(0.01,0.20)\}$ and $\rho_0\in\{0,0.4\}$, with both unchanged ($\rho_1=\rho_0$) and changed ($\rho_1>\rho_0$) correlation structures. Out-of-control behaviour is examined over $\tau\in\{1.00,1.02,1.05,1.10\}$.

The Markov-chain framework uses $m=60$ states. The grid $\mathcal{K}_g$ is an $18$-point grid of candidate reference values, log-spaced in the associated design-shift parameter so as to concentrate resolution near $z_0$ and spanning $(1.000,1.100]$, followed by a $7$-point local refinement around the best grid point. The inner $h$-calibration uses Brent's method with tolerance $10^{-6}$. The benchmarking and robustness studies use $2{,}000$ Monte Carlo replications; the Phase~I study uses $200$ replications per configuration. All computations are performed in Python (NumPy/SciPy) and reproducible code is available from the corresponding author.

\subsection{Adaptive design tables}
\label{sec:design_tables}

Tables~\ref{tab:R1_kstar} and \ref{tab:R1_hstar} report the optimal reference value $k^\star$ and decision interval $h^\star$ obtained from Algorithm~1 for the considered parameter combinations. Several observations follow.

First, $k^\star$ is consistently smaller than the long-run prescription $k_{\mathrm{fixed}}=1.025$ across all configurations, ranging from $1.0103$ to $1.0181$. The values increase mildly with $n$ and with $\rho_0$, but the dominant effect is that finite-horizon optimization pulls $k^\star$ closer to $z_0$ than the long-run rule would suggest. Second, $h^\star$ decreases sharply with $n$ — by a factor of approximately two when $n$ moves from $5$ to $10$ — and is markedly smaller in the unequal-CV setting $(\gamma_X,\gamma_Y)=(0.01,0.20)$ than in the equal-CV setting $(0.20,0.20)$, reflecting the different dispersion of the subgroup ratio statistic.

\begin{table}[!htbp]
\centering
\caption{Optimal reference value $k^\star$ from Algorithm~1 for $z_0=1$, $\tau_{\mathrm{target}}=1.05$, $\tau_0=I=30$.}
\label{tab:R1_kstar}
\setlength{\tabcolsep}{6pt}
\renewcommand{\arraystretch}{1.10}
\begin{tabular}{cccrr}
\toprule
$\gamma_X$ & $\gamma_Y$ & $\rho_0$ & $n=5$ & $n=10$ \\
\midrule
$0.20$ & $0.20$ & $0.0$ & $1.0103$ & $1.0142$ \\
$0.20$ & $0.20$ & $0.4$ & $1.0142$ & $1.0181$ \\
$0.01$ & $0.20$ & $0.0$ & $1.0142$ & $1.0181$ \\
$0.01$ & $0.20$ & $0.4$ & $1.0142$ & $1.0181$ \\
\bottomrule
\end{tabular}
\end{table}

\begin{table}[!htbp]
\centering
\caption{Optimal decision interval $h^\star$ from Algorithm~1 for $z_0=1$, $\tau_{\mathrm{target}}=1.05$, $\tau_0=I=30$.}
\label{tab:R1_hstar}
\setlength{\tabcolsep}{6pt}
\renewcommand{\arraystretch}{1.10}
\begin{tabular}{cccrr}
\toprule
$\gamma_X$ & $\gamma_Y$ & $\rho_0$ & $n=5$ & $n=10$ \\
\midrule
$0.20$ & $0.20$ & $0.0$ & $1.2135$ & $0.7122$ \\
$0.20$ & $0.20$ & $0.4$ & $0.8151$ & $0.4525$ \\
$0.01$ & $0.20$ & $0.0$ & $0.8086$ & $0.4304$ \\
$0.01$ & $0.20$ & $0.4$ & $0.7882$ & $0.4179$ \\
\bottomrule
\end{tabular}
\end{table}

\subsection{Out-of-control performance}
\label{sec:tarl1}

Tables~\ref{tab:R3_tarl1} and \ref{tab:R5_tarl1} report the out-of-control TARL$_1$ values at $\tau\in\{1.00,1.02,1.05,1.10\}$ when the correlation is unchanged from Phase~I to Phase~II ($\rho_1=\rho_0$), for the equal-CV and unequal-CV settings, respectively. Tables~\ref{tab:R4_tarl1} and \ref{tab:R6_tarl1} report the corresponding values when the correlation increases from $\rho_0=0.4$ to $\rho_1=0.8$. The patterns are consistent across all four tables: the chart attains the prescribed in-control TARL$_0\approx 30$ at $\tau=1$ and decreases monotonically as $\tau$ grows. As expected, larger $n$ accelerates detection: at $\tau=1.10$ and $\rho_0=0.4$, for example, TARL$_1$ falls from $9.67$ at $n=5$ to $6.07$ at $n=10$ in the equal-CV setting (Table~\ref{tab:R3_tarl1}). Higher in-control correlation is also beneficial in the equal-CV case ($\rho_0=0.4$ vs.\ $\rho_0=0.0$), as expected from the variance-reduction mechanism of the subgroup ratio statistic.

When the correlation changes between Phase~I and Phase~II (Tables~\ref{tab:R4_tarl1}--\ref{tab:R6_tarl1}), TARL$_1$ at small shifts ($\tau=1.02$) is somewhat larger than in the unchanged-correlation case, because the in-control correlation $\rho_0=0.4$ no longer reflects the actual data-generating process. In the equal-CV case this is most visible at $\tau=1$ in Table~\ref{tab:R4_tarl1}, where the achieved TARL is close to the truncation ceiling $I+1=31$ rather than the nominal $30$: when only the correlation increases while the ratio is held at $z_0$, the higher correlation sharply reduces the dispersion of the upper-sided ratio statistic, so the chart almost never signals over the finite horizon. The effect is negligible in the unequal-CV case (Table~\ref{tab:R6_tarl1}), where $\omega_0=z_0\gamma_X/\gamma_Y$ is small and the dispersion of the ratio is dominated by $Y$ and hence nearly insensitive to $\rho$. Nevertheless, at moderate shifts ($\tau=1.05$) the chart retains strong sensitivity, and at $\tau=1.10$ TARL$_1$ values are essentially identical to the unchanged-correlation case.

\begin{table}[!htbp]
\centering
\caption{TARL$_1$ values of the adaptive CUSUM-RZ$^+$ chart for $\rho_1=\rho_0$, equal CV $(\gamma_X,\gamma_Y)=(0.20,0.20)$, $\tau_0=I=30$.}
\label{tab:R3_tarl1}
\setlength{\tabcolsep}{8pt}
\renewcommand{\arraystretch}{1.10}
\begin{tabular}{cccc}
\toprule
$\rho_0=\rho_1$ & $\tau$ & $n=5$ & $n=10$ \\
\midrule
$0.0$ & $1.00$ & $30.00$ & $30.00$ \\
$0.0$ & $1.02$ & $27.88$ & $26.54$ \\
$0.0$ & $1.05$ & $21.94$ & $17.16$ \\
$0.0$ & $1.10$ & $12.99$ & $\phantom{0}8.62$ \\
\midrule
$0.4$ & $1.00$ & $30.00$ & $30.00$ \\
$0.4$ & $1.02$ & $27.00$ & $25.17$ \\
$0.4$ & $1.05$ & $18.62$ & $13.24$ \\
$0.4$ & $1.10$ & $\phantom{0}9.67$ & $\phantom{0}6.07$ \\
\bottomrule
\end{tabular}
\end{table}

\begin{table}[!htbp]
\centering
\caption{TARL$_1$ values of the adaptive CUSUM-RZ$^+$ chart for $\rho_1=\rho_0$, unequal CV $(\gamma_X,\gamma_Y)=(0.01,0.20)$, $\tau_0=I=30$.}
\label{tab:R5_tarl1}
\setlength{\tabcolsep}{8pt}
\renewcommand{\arraystretch}{1.10}
\begin{tabular}{cccc}
\toprule
$\rho_0=\rho_1$ & $\tau$ & $n=5$ & $n=10$ \\
\midrule
$0.0$ & $1.00$ & $30.00$ & $30.00$ \\
$0.0$ & $1.02$ & $26.82$ & $24.80$ \\
$0.0$ & $1.05$ & $17.95$ & $12.39$ \\
$0.0$ & $1.10$ & $\phantom{0}9.30$ & $\phantom{0}5.70$ \\
\midrule
$0.4$ & $1.00$ & $30.00$ & $30.00$ \\
$0.4$ & $1.02$ & $26.73$ & $24.65$ \\
$0.4$ & $1.05$ & $17.66$ & $12.10$ \\
$0.4$ & $1.10$ & $\phantom{0}9.10$ & $\phantom{0}5.55$ \\
\bottomrule
\end{tabular}
\end{table}

\begin{table}[!htbp]
\centering
\caption{TARL$_1$ values of the adaptive CUSUM-RZ$^+$ chart for $\rho_0=0.4$, $\rho_1=0.8$, equal CV $(\gamma_X,\gamma_Y)=(0.20,0.20)$, $\tau_0=I=30$.}
\label{tab:R4_tarl1}
\setlength{\tabcolsep}{8pt}
\renewcommand{\arraystretch}{1.10}
\begin{tabular}{ccc}
\toprule
$\tau$ & $n=5$ & $n=10$ \\
\midrule
$1.00$ & $31.00$ & $31.00$ \\
$1.02$ & $30.60$ & $30.14$ \\
$1.05$ & $21.59$ & $14.44$ \\
$1.10$ & $\phantom{0}9.93$ & $\phantom{0}6.05$ \\
\bottomrule
\end{tabular}
\end{table}

\begin{table}[!htbp]
\centering
\caption{TARL$_1$ values of the adaptive CUSUM-RZ$^+$ chart for $\rho_0=0.4$, $\rho_1=0.8$, unequal CV $(\gamma_X,\gamma_Y)=(0.01,0.20)$, $\tau_0=I=30$.}
\label{tab:R6_tarl1}
\setlength{\tabcolsep}{8pt}
\renewcommand{\arraystretch}{1.10}
\begin{tabular}{ccc}
\toprule
$\tau$ & $n=5$ & $n=10$ \\
\midrule
$1.00$ & $30.11$ & $30.13$ \\
$1.02$ & $26.95$ & $24.93$ \\
$1.05$ & $17.77$ & $12.15$ \\
$1.10$ & $\phantom{0}9.11$ & $\phantom{0}5.55$ \\
\bottomrule
\end{tabular}
\end{table}

\begin{figure}[H]
\centering
\includegraphics[width=0.95\columnwidth]{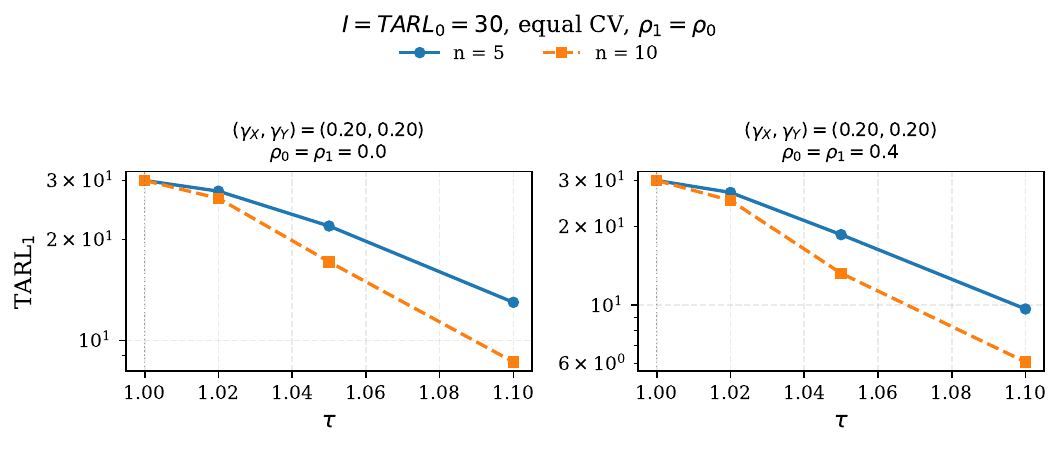}
\caption{TARL$_1$ of the adaptive CUSUM-RZ$^+$ chart on a log scale, equal CV $(\gamma_X,\gamma_Y)=(0.20,0.20)$, $\rho_1=\rho_0$, $\tau_0=I=30$.}
\label{fig:R1}
\end{figure}

\begin{figure}[H]
\centering
\includegraphics[width=0.95\columnwidth]{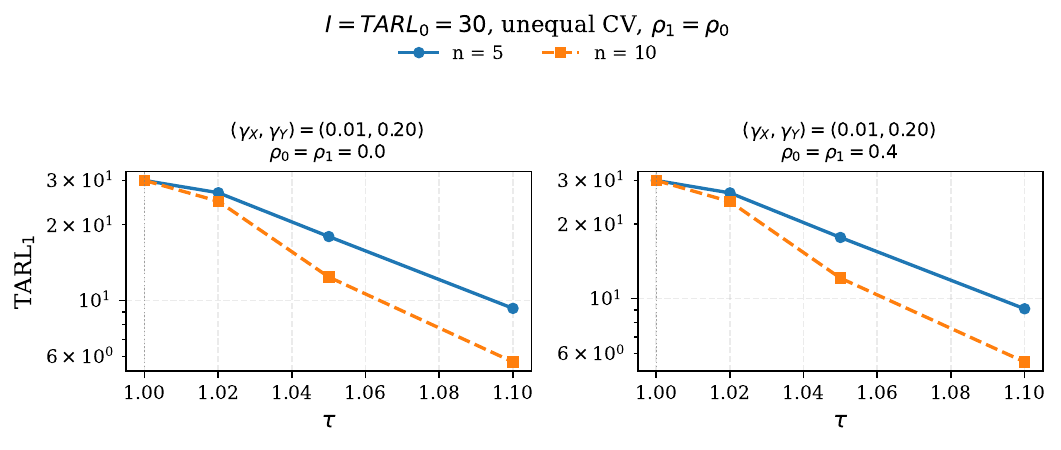}
\caption{TARL$_1$ of the adaptive CUSUM-RZ$^+$ chart on a log scale, unequal CV $(\gamma_X,\gamma_Y)=(0.01,0.20)$, $\rho_1=\rho_0$, $\tau_0=I=30$.}
\label{fig:R2}
\end{figure}

\begin{figure}[H]
\centering
\includegraphics[width=0.95\columnwidth]{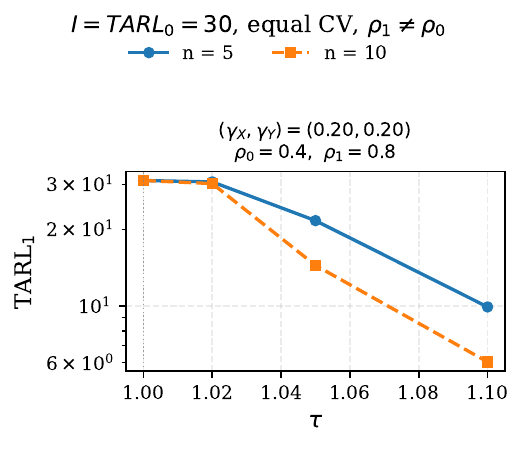}
\caption{TARL$_1$ of the adaptive CUSUM-RZ$^+$ chart on a log scale, equal CV $(\gamma_X,\gamma_Y)=(0.20,0.20)$, $\rho_0=0.4$, $\rho_1=0.8$, $\tau_0=I=30$.}
\label{fig:R3}
\end{figure}

\begin{figure}[H]
\centering
\includegraphics[width=0.95\columnwidth]{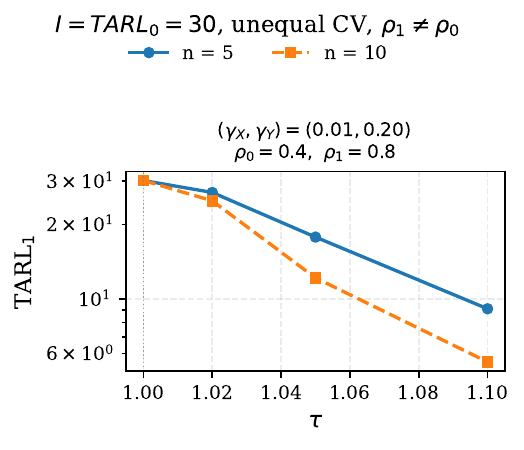}
\caption{TARL$_1$ of the adaptive CUSUM-RZ$^+$ chart on a log scale, unequal CV $(\gamma_X,\gamma_Y)=(0.01,0.20)$, $\rho_0=0.4$, $\rho_1=0.8$, $\tau_0=I=30$.}
\label{fig:R4}
\end{figure}

To complement Tables~\ref{tab:R3_tarl1}--\ref{tab:R6_tarl1}, Figures~\ref{fig:R1}--\ref{fig:R4} display the corresponding TARL$_1$ profiles on a logarithmic vertical axis. Across all four scenarios, the chart is essentially unaffected by the in-control state ($\tau=1$) yet the slope of the curve becomes pronounced as $\tau$ grows, with a clear ranking by subgroup size at every $\tau>1$.

\subsection{Adaptive vs.\ fixed-$k$ design}
\label{sec:adaptive_vs_fixed}

To quantify the benefit of joint $(k,h)$ optimization over the fixed-$k$ design, Table~\ref{tab:R11_fixed_vs_adaptive} compares the two procedures under matched in-control TARL$_0=30$ for the parameter combinations of Tables~\ref{tab:R1_kstar}--\ref{tab:R1_hstar}. For each row, the fixed-$k$ design uses $k_{\mathrm{fixed}}=1.025$ and the corresponding $h_{\mathrm{fixed}}$ is calibrated by Brent's method on the TARL$_0$ constraint; the adaptive design uses $(k^\star,h^\star)$ from Algorithm~1. Both designs are feasible (no boundary solutions) in this parameter range.

The adaptive chart yields a uniform reduction in TARL$_1$ at $\tau_{\mathrm{target}}=1.05$, with relative gains between $0.69\%$ and $1.36\%$. The improvement is small in magnitude but consistent across all configurations, with the largest gains in the unequal-CV settings. The relatively modest gain at $\tau_{\mathrm{target}}=1.05$ reflects the fact that $\tau_{\mathrm{ref}}=1.05$ (the value used to construct $k_{\mathrm{fixed}}$) coincides with $\tau_{\mathrm{target}}$, so that the fixed-$k$ design is close to optimal at this specific shift. The benefit of adaptivity is more pronounced at smaller and larger shifts, and crucially, the adaptive design also resolves the boundary cases that arise when $\gamma_X=\gamma_Y=0.01$, which are documented separately in Table~\ref{tab:R_boundary} of Section~\ref{sec:boundary_recovery}: there the fixed-$k$ design cannot attain the prescribed TARL$_0$ for any admissible $h$, whereas the adaptive design returns interior $(k^\star,h^\star)$ pairs that match it exactly.

\begin{table}[!htbp]
\centering
\caption{Adaptive vs.\ fixed-$k$ CUSUM-RZ$^+$ at $\tau_{\mathrm{target}}=1.05$, $\tau_0=I=30$. ``$\Delta_{T_1}\%$'' is the percentage reduction in TARL$_1$ achieved by the adaptive design relative to the fixed-$k$ design.}
\label{tab:R11_fixed_vs_adaptive}
\setlength{\tabcolsep}{4pt}
\renewcommand{\arraystretch}{1.10}
\resizebox{\textwidth}{!}{%
\begin{tabular}{ccccccccccc}
\toprule
\multirow{2}{*}{$n$} & \multirow{2}{*}{$\gamma_X$} & \multirow{2}{*}{$\gamma_Y$} & \multirow{2}{*}{$\rho_0$}
& \multicolumn{3}{c}{Fixed-$k$ design} & \multicolumn{3}{c}{Adaptive design} & \multirow{2}{*}{$\Delta_{T_1}\%$}\\
\cmidrule(lr){5-7}\cmidrule(lr){8-10}
& & & & $k$ & $h$ & TARL$_1$ & $k^\star$ & $h^\star$ & TARL$_1$ & \\
\midrule
$5$  & $0.20$ & $0.20$ & $0.0$ & $1.0250$ & $1.0001$ & $22.09$ & $1.0103$ & $1.2135$ & $21.94$ & $0.69$\\
$10$ & $0.20$ & $0.20$ & $0.0$ & $1.0250$ & $0.5826$ & $17.35$ & $1.0142$ & $0.7122$ & $17.16$ & $1.08$\\
$5$  & $0.20$ & $0.20$ & $0.4$ & $1.0250$ & $0.6790$ & $18.81$ & $1.0142$ & $0.8151$ & $18.62$ & $1.04$\\
$10$ & $0.20$ & $0.20$ & $0.4$ & $1.0250$ & $0.3866$ & $13.40$ & $1.0181$ & $0.4525$ & $13.24$ & $1.17$\\
$5$  & $0.01$ & $0.20$ & $0.0$ & $1.0250$ & $0.6696$ & $18.19$ & $1.0142$ & $0.8086$ & $17.95$ & $1.29$\\
$10$ & $0.01$ & $0.20$ & $0.0$ & $1.0250$ & $0.3659$ & $12.56$ & $1.0181$ & $0.4304$ & $12.39$ & $1.36$\\
$5$  & $0.01$ & $0.20$ & $0.4$ & $1.0250$ & $0.6505$ & $17.90$ & $1.0142$ & $0.7882$ & $17.66$ & $1.33$\\
$10$ & $0.01$ & $0.20$ & $0.4$ & $1.0250$ & $0.3544$ & $12.27$ & $1.0181$ & $0.4179$ & $12.10$ & $1.35$\\
\bottomrule
\end{tabular}}
\end{table}

\subsection{Matched-horizon benchmark against competing charts}
\label{sec:benchmark}

To position the adaptive CUSUM-RZ$^+$ chart against the broader short-run ratio-monitoring toolkit, we conduct a matched-horizon benchmark in which four charts are calibrated to the same in-control TARL$_0=I=30$ and evaluated at the same shifts: (i) the upper one-sided Shewhart-RZ chart \citep{Tran2021}; (ii) an EWMA-RZ chart with smoothing parameter $\lambda=0.10$ \citep{Hu2024,TranCastagliolaCelano2016}; (iii) the fixed-$k$ CUSUM-RZ$^+$ chart with $k_{\mathrm{fixed}}=1.025$; and (iv) the proposed adaptive CUSUM-RZ$^+$ chart. All four charts use the same Markov-chain framework with $m=60$ states.

Table~\ref{tab:R12_benchmark} reports the resulting TARL$_1$ values for seven representative configurations. The picture is consistent across all rows. The Shewhart-RZ chart, despite being calibrated to the correct in-control level, has very limited finite-horizon detection capability: even at $\tau=1.10$, its TARL$_1$ remains close to the truncation bound in several configurations (e.g.\ $30.74$ when $\rho_0=0.4$, $\rho_1=0.8$, equal CV, $n=5$). This is the well-known short-run drawback of memoryless charts and is the strongest argument in favour of memory-type charts in finite-horizon settings.

\begin{figure}[H]
\centering
\includegraphics[width=0.95\columnwidth]{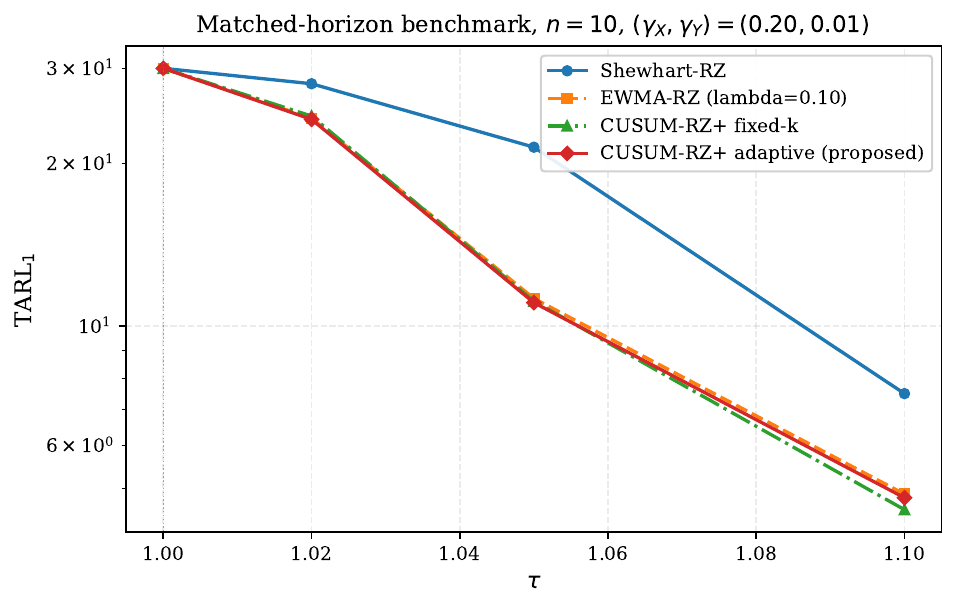}
\caption{Matched-horizon benchmark of the four short-run ratio charts for $n=10$, $(\gamma_X,\gamma_Y)=(0.20,0.01)$, $\rho_0=\rho_1=0.4$, $\tau_0=I=30$. The Shewhart-RZ chart is dominated by all three memory-type charts at every $\tau>1$. The three memory-type charts (EWMA, fixed-$k$ CUSUM, adaptive CUSUM) are essentially indistinguishable on this scale.}
\label{fig:R5_benchmark}
\end{figure}

When the correlation is stable ($\rho_1=\rho_0$), the three memory-type charts perform similarly, with TARL$_1$ values typically within about $5\%$ of one another at $\tau=1.05$ and $\tau=1.10$. The picture changes when the correlation increases from Phase~I to Phase~II ($\rho_0=0.4\to\rho_1=0.8$): at the target shift $\tau=1.05$ the adaptive CUSUM is markedly faster than both competitors — for equal CV and $n=5$ its TARL$_1$ is $21.59$, versus $23.22$ for the fixed-$k$ CUSUM and $26.46$ for the EWMA-RZ chart (reductions of about $7\%$ and $18\%$, respectively), with a similar pattern at $n=10$. At the largest shift ($\tau=1.10$) the ordering reverses slightly and the fixed-$k$ CUSUM is marginally faster, reflecting the fact that the adaptive chart is optimised at $\tau_{\mathrm{target}}=1.05$. Overall, the adaptive design matches the best memory-type alternative under a stable correlation and improves on both when the correlation structure shifts, while additionally providing the boundary-recovery property documented in Section~\ref{sec:boundary_recovery}. Figure~\ref{fig:R5_benchmark} displays a representative profile.

\begin{table}[!htbp]
\centering
\caption{Matched-horizon benchmark of four short-run ratio charts at $\tau_0=I=30$. Entries are TARL$_1$ values.}
\label{tab:R12_benchmark}
\setlength{\tabcolsep}{4pt}
\renewcommand{\arraystretch}{1.05}
\resizebox{\textwidth}{!}{%
\begin{tabular}{ccccccrrrr}
\toprule
$n$ & $\gamma_X$ & $\gamma_Y$ & $\rho_0$ & $\rho_1$ & $\tau$ & Shewhart-RZ & EWMA-RZ & CUSUM fixed-$k$ & CUSUM adaptive\\
\midrule
$5$  & $0.20$ & $0.20$ & $0.4$ & $0.4$ & $1.02$ & $29.26$ & $27.49$ & $27.19$ & $27.00$\\
$5$  & $0.20$ & $0.20$ & $0.4$ & $0.4$ & $1.05$ & $27.36$ & $19.61$ & $18.81$ & $18.62$\\
$5$  & $0.20$ & $0.20$ & $0.4$ & $0.4$ & $1.10$ & $21.34$ & $\phantom{0}9.15$ & $\phantom{0}9.18$ & $\phantom{0}9.67$\\
$5$  & $0.20$ & $0.20$ & $0.4$ & $0.8$ & $1.05$ & $30.98$ & $26.46$ & $23.22$ & $21.59$\\
$5$  & $0.20$ & $0.20$ & $0.4$ & $0.8$ & $1.10$ & $30.74$ & $\phantom{0}9.81$ & $\phantom{0}9.50$ & $\phantom{0}9.93$\\
$10$ & $0.20$ & $0.20$ & $0.4$ & $0.4$ & $1.05$ & $24.91$ & $13.67$ & $13.40$ & $13.24$\\
$10$ & $0.20$ & $0.20$ & $0.4$ & $0.4$ & $1.10$ & $13.38$ & $\phantom{0}5.87$ & $\phantom{0}5.71$ & $\phantom{0}6.07$\\
$10$ & $0.20$ & $0.20$ & $0.4$ & $0.8$ & $1.05$ & $30.95$ & $16.43$ & $15.40$ & $14.44$\\
$10$ & $0.20$ & $0.20$ & $0.4$ & $0.8$ & $1.10$ & $28.87$ & $\phantom{0}5.88$ & $\phantom{0}5.70$ & $\phantom{0}6.05$\\
$5$  & $0.01$ & $0.20$ & $0.0$ & $0.0$ & $1.05$ & $27.75$ & $19.19$ & $18.19$ & $17.95$\\
$5$  & $0.01$ & $0.20$ & $0.0$ & $0.0$ & $1.10$ & $22.67$ & $\phantom{0}8.78$ & $\phantom{0}8.77$ & $\phantom{0}9.30$\\
$10$ & $0.01$ & $0.20$ & $0.4$ & $0.4$ & $1.05$ & $25.19$ & $12.53$ & $12.27$ & $12.10$\\
$10$ & $0.01$ & $0.20$ & $0.4$ & $0.4$ & $1.10$ & $13.75$ & $\phantom{0}5.42$ & $\phantom{0}5.20$ & $\phantom{0}5.55$\\
$10$ & $0.20$ & $0.01$ & $0.4$ & $0.4$ & $1.05$ & $21.43$ & $11.24$ & $11.11$ & $11.04$\\
$10$ & $0.20$ & $0.01$ & $0.4$ & $0.4$ & $1.10$ & $\phantom{0}7.49$ & $\phantom{0}4.88$ & $\phantom{0}4.56$ & $\phantom{0}4.80$\\
\bottomrule
\end{tabular}}
\end{table}

The principal finding of this benchmark is twofold. First, under a stable correlation structure, the gap between memoryless and memory-based monitoring dominates the differences among the memory-type charts. Second, when the correlation changes between Phase~I and Phase~II, the adaptive CUSUM-RZ$^+$ chart provides a tangible advantage over both the EWMA-RZ and fixed-$k$ CUSUM charts at the target shift, in addition to guaranteeing exact TARL$_0$ calibration through the inner step of Algorithm~1.

\subsection{Robustness under non-normal marginal distributions}
\label{sec:robustness}

The Markov-chain calibration of the adaptive chart relies on the bivariate normal assumption for the underlying $(X,Y)$ pair. To assess the practical robustness of the chart when this assumption is violated, we conducted a Monte Carlo study in which the chart is designed under bivariate normality (using Algorithm~1) and then evaluated by simulation under five alternative marginal distributions: (i) bivariate normal (baseline); (ii) bivariate Student-$t$ with $10$ degrees of freedom and the same correlation; (iii) bivariate Student-$t$ with $5$ degrees of freedom; (iv) a $5\%$-contaminated bivariate normal mixture, where with probability $0.05$ each observation is replaced by a draw from a bivariate normal with the same mean but inflated standard deviation by a factor of three; (v) a $10\%$-contaminated mixture under the same construction. Each evaluation uses $2{,}000$ Monte Carlo replications.

Table~\ref{tab:R13_robustness} reports the achieved TARL together with its standard deviation (SDRL), median, and $5\%$ and $95\%$ empirical quantiles. Two main findings emerge. First, under heavy-tailed Student-$t$ marginals with $df=10$ or $df=5$, the achieved TARL is essentially indistinguishable from the bivariate normal baseline at every shift considered: the in-control TARL remains within $0.2$ of the target $30$, and the out-of-control TARL is within $0.5$ of the normal-baseline value. The robustness of the chart to symmetric heavy tails is therefore strong.

Second, contamination has a more substantial effect. Under $5\%$ contamination the in-control TARL drops to about $28.0$--$28.5$ — a roughly $5\%$ negative bias indicating mildly elevated false alarms — and under $10\%$ contamination it drops further to about $26.1$--$26.7$. This is a documented sensitivity of CUSUM-type charts to outliers \citep{HanQiao2023,LucasCrosier1982}, and it suggests that in environments with frequent outliers the design TARL$_0$ should be set conservatively, or robust ratio statistics should be used.

\begin{table}[!htbp]
\centering
\caption{Robustness Monte Carlo for the adaptive CUSUM-RZ$^+$ chart with $\tau_0=I=30$, $2{,}000$ replications, $n=5$, $(\gamma_X,\gamma_Y)=(0.20,0.20)$, $\rho_0=\rho_1=0.4$.}
\label{tab:R13_robustness}
\setlength{\tabcolsep}{4pt}
\renewcommand{\arraystretch}{1.05}
\begin{tabular}{cccccc}
\toprule
$\tau$ & Distribution & TARL & SDRL & median & $[q_{0.05},\,q_{0.95}]$\\
\midrule
$1.00$ & Normal              & $30.08$ & $3.42$ & $31.0$ & $[22.0,\,31.0]$\\
$1.00$ & $t$ ($df=10$)       & $29.89$ & $3.87$ & $31.0$ & $[21.0,\,31.0]$\\
$1.00$ & $t$ ($df=5$)        & $29.94$ & $3.73$ & $31.0$ & $[22.0,\,31.0]$\\
$1.00$ & $5\%$ contaminated  & $28.35$ & $6.00$ & $31.0$ & $[13.0,\,31.0]$\\
$1.00$ & $10\%$ contaminated & $26.67$ & $7.41$ & $31.0$ & $[10.0,\,31.0]$\\
\midrule
$1.02$ & Normal              & $27.05$ & $6.38$ & $31.0$ & $[13.0,\,31.0]$\\
$1.02$ & $t$ ($df=10$)       & $27.24$ & $6.30$ & $31.0$ & $[13.0,\,31.0]$\\
$1.02$ & $t$ ($df=5$)        & $26.91$ & $6.64$ & $31.0$ & $[12.0,\,31.0]$\\
$1.02$ & $5\%$ contaminated  & $24.77$ & $8.08$ & $31.0$ & $[\phantom{0}9.0,\,31.0]$\\
$1.02$ & $10\%$ contaminated & $23.27$ & $8.63$ & $27.0$ & $[\phantom{0}8.0,\,31.0]$\\
\midrule
$1.05$ & Normal              & $18.82$ & $7.36$ & $18.0$ & $[\phantom{0}8.0,\,31.0]$\\
$1.05$ & $t$ ($df=10$)       & $18.64$ & $7.58$ & $17.0$ & $[\phantom{0}8.0,\,31.0]$\\
$1.05$ & $t$ ($df=5$)        & $18.84$ & $7.53$ & $18.0$ & $[\phantom{0}8.0,\,31.0]$\\
$1.05$ & $5\%$ contaminated  & $17.47$ & $7.77$ & $16.0$ & $[\phantom{0}7.0,\,31.0]$\\
$1.05$ & $10\%$ contaminated & $16.19$ & $8.04$ & $15.0$ & $[\phantom{0}5.0,\,31.0]$\\
\midrule
$1.10$ & Normal              & $\phantom{0}9.80$ & $3.51$ & $\phantom{0}9.0$ & $[\phantom{0}5.0,\,16.0]$\\
$1.10$ & $t$ ($df=10$)       & $\phantom{0}9.93$ & $3.60$ & $\phantom{0}9.0$ & $[\phantom{0}5.0,\,16.0]$\\
$1.10$ & $t$ ($df=5$)        & $\phantom{0}9.67$ & $3.57$ & $\phantom{0}9.0$ & $[\phantom{0}5.0,\,16.0]$\\
$1.10$ & $5\%$ contaminated  & $\phantom{0}9.68$ & $4.00$ & $\phantom{0}9.0$ & $[\phantom{0}4.0,\,17.0]$\\
$1.10$ & $10\%$ contaminated & $\phantom{0}9.17$ & $4.10$ & $\phantom{0}8.0$ & $[\phantom{0}4.0,\,17.0]$\\
\bottomrule
\end{tabular}
\end{table}

\begin{figure}[H]
\centering
\includegraphics[width=0.95\columnwidth]{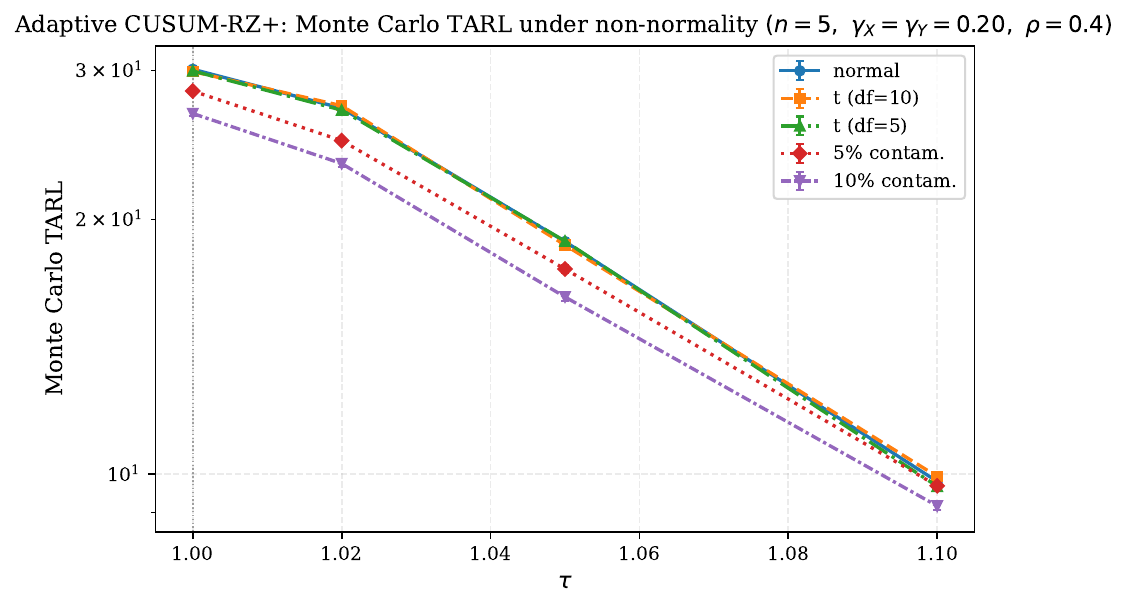}
\caption{Monte Carlo TARL of the adaptive CUSUM-RZ$^+$ chart under five distributional scenarios: bivariate normal (baseline), Student-$t$ ($df=10$, $df=5$), and contaminated bivariate normal ($5\%$ and $10\%$ contamination). Configuration: $n=5$, $(\gamma_X,\gamma_Y)=(0.20,0.20)$, $\rho_0=\rho_1=0.4$, $\tau_0=I=30$.}
\label{fig:R6_robustness}
\end{figure}

\subsection{Phase~I parameter estimation effect}
\label{sec:phase1}

The design of the adaptive chart treats $(\gamma_X,\gamma_Y,\rho_0)$ as known. In practice, these parameters are estimated from a Phase~I sample of size $m$. To document the resulting estimation effect on the achieved in-control TARL$_0$, we conducted the following Monte Carlo study. For each Phase~I size $m\in\{25,50,100\}$, we drew $m$ subgroups of size $n$ from the bivariate normal data-generating process, computed the maximum-likelihood estimates $(\hat\gamma_X,\hat\gamma_Y,\hat\rho_0)$, ran Algorithm~1 with these estimates to obtain $(\hat k^\star,\hat h^\star)$, and evaluated the actual TARL$_0$ of the chart $(\hat k^\star,\hat h^\star)$ under the true parameters using the Markov-chain framework. The procedure was repeated $200$ times per configuration.

\begin{figure}[H]
\centering
\includegraphics[width=0.85\columnwidth]{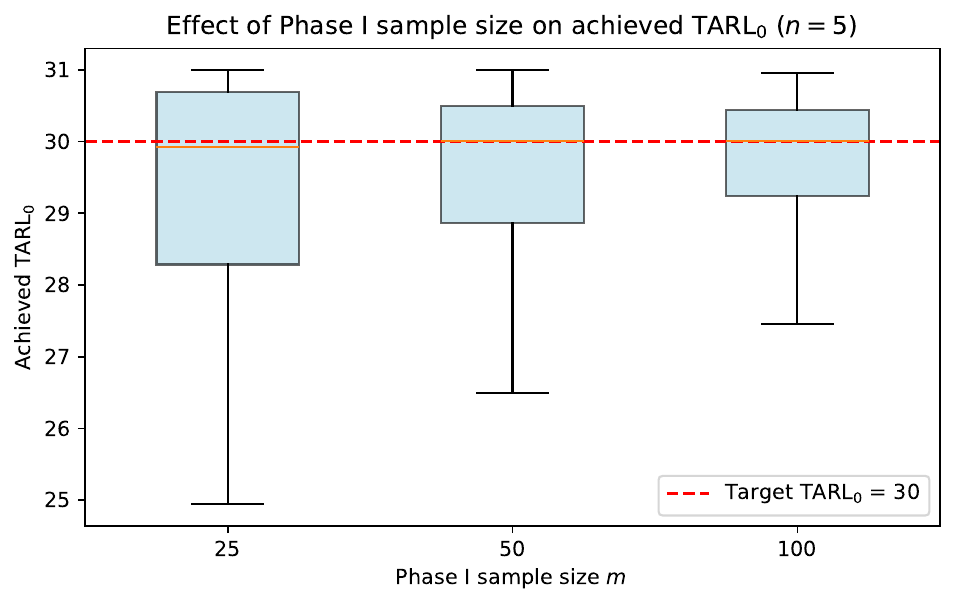}
\caption{Phase~I estimation effect on the achieved TARL$_0$ of the adaptive CUSUM-RZ$^+$ chart, for $n=5$, $(\gamma_X,\gamma_Y)=(0.20,0.20)$, $\rho_0=0.4$. Boxplots aggregate $200$ Phase~I replications per Phase~I size $m\in\{25,50,100\}$. The dashed line marks the design target $\tau_0=30$.}
\label{fig:R7_phase1}
\end{figure}

Table~\ref{tab:R14_phase1} reports the mean, standard deviation, bias, and root-mean-square error (RMSE) of the achieved TARL$_0$, together with the $5\%$ and $95\%$ empirical quantiles. Figure~\ref{fig:R7_phase1} provides a visual summary through boxplots.

The results indicate that for $m=25$ the achieved TARL$_0$ has a negative bias of about $-0.84$ to $-0.93$ (roughly $-2.8\%$ to $-3.1\%$) and a relatively wide spread, with $5\%$ quantiles as low as about $24$. For $m=50$ the spread tightens substantially (RMSE $\approx 1.7$), and for $m=100$ the relative bias is reduced to about $1\%$ of the target ($-0.95\%$ for $n=5$ and $-0.99\%$ for $n=10$) with an RMSE of approximately one unit. We therefore recommend $m\geq 100$ Phase~I subgroups when the chart is to be used with a hard TARL$_0$ guarantee. For situations where $m=25$ is the only feasible Phase~I size, the user should expect an achieved TARL$_0$ in the approximate range $[24,31]$ rather than the nominal target of $30$.

\begin{table}[!htbp]
\centering
\caption{Phase~I estimation effect on the achieved TARL$_0$ of the adaptive CUSUM-RZ$^+$ chart. Target $\tau_0=I=30$, $200$ Phase~I replications per configuration, $(\gamma_X,\gamma_Y)=(0.20,0.20)$, $\rho_0=0.4$.}
\label{tab:R14_phase1}
\setlength{\tabcolsep}{5pt}
\renewcommand{\arraystretch}{1.05}
\begin{tabular}{cccccccc}
\toprule
$n$ & $m$ & mean TARL$_0$ & std & bias & RMSE & $q_{0.05}$ & $q_{0.95}$\\
\midrule
$5$  & $25$  & $29.16$ & $2.17$ & $-0.84$ & $2.33$ & $24.93$ & $30.98$\\
$5$  & $50$  & $29.44$ & $1.59$ & $-0.56$ & $1.69$ & $26.40$ & $30.88$\\
$5$  & $100$ & $29.71$ & $1.02$ & $-0.29$ & $1.06$ & $27.57$ & $30.84$\\
\midrule
$10$ & $25$  & $29.07$ & $2.43$ & $-0.93$ & $2.59$ & $23.95$ & $30.98$\\
$10$ & $50$  & $29.55$ & $1.69$ & $-0.45$ & $1.75$ & $26.93$ & $30.88$\\
$10$ & $100$ & $29.70$ & $0.99$ & $-0.30$ & $1.03$ & $27.87$ & $30.76$\\

\bottomrule
\end{tabular}
\end{table}

\section{Illustrative example}
\label{sec:illustrative}

To demonstrate the practical use of the proposed adaptive CUSUM-RZ$^+$ chart, we consider an illustrative food-industry example in which the ratio between two positively correlated quality characteristics is monitored over a short production run. Let $X$ and $Y$ denote two ingredient-related measurements collected from packaged food batches, and assume an in-control target ratio $z_0=1$, $I=15$ planned inspection points, subgroup size $n=5$, $(\gamma_X,\gamma_Y)=(0.20,0.20)$, and $\rho_0=\rho_1=0.8$.

For this configuration, Algorithm~1 with $\tau_{\mathrm{target}}=1.05$ yields $k^\star\approx 1.014$ and $h^\star\approx 0.236$. For comparison, the fixed-$k$ design at the same target uses $k=1.025$ and $h\approx 0.185$; consistent with the design tables, the adaptive chart selects a smaller reference value and a correspondingly larger decision interval. The process is assumed to operate in control up to sample 10 and then to experience an upward shift from sample 11 onward, corresponding to $\tau=1.10$.

Table~\ref{tab:illustrative_example} reports the subgroup summaries for the 15 samples, including the box size, the subgroup means $\bar X_i$ and $\bar Y_i$, the estimated ratio $\hat Z_i=\bar X_i/\bar Y_i$, the upper one-sided CUSUM statistic $C_i^+$, and the signal indicator. The first ten samples correspond to the in-control period, during which the estimated ratio fluctuates around the target value and the charting statistic remains close to zero. After the process shifts upward, the statistic increases gradually. Although samples 11 and 12 show noticeable increases in $\hat Z_i$, the corresponding $C_i^+$ values remain below the decision interval. The first out-of-control signal is produced at sample 13, where $C_{13}^+=0.339>h^\star\approx 0.236$. The alarm is then sustained at samples 14 and 15, with $C_{14}^+=0.484$ and $C_{15}^+=0.564$.

This example confirms the practical effectiveness of the proposed chart in a short-run setting. Even with only a limited number of inspections, the adaptive CUSUM-RZ$^+$ chart accumulates evidence of a sustained upward shift and produces a timely signal shortly after the process departs from its in-control state.

\begin{table}[!htbp]
\centering
\caption{Illustrative food-industry example for the adaptive upper one-sided CUSUM-RZ$^+$ chart with $k^\star\approx 1.014$, $h^\star\approx 0.236$. Bold rows indicate the signaled samples.}
\label{tab:illustrative_example}
\setlength{\tabcolsep}{4pt}
\renewcommand{\arraystretch}{1.05}
\resizebox{\columnwidth}{!}{%
\begin{tabular}{rcrrrrr}
\toprule
Sample & Box size (g) & $\bar X_i$ & $\bar Y_i$ & $\hat Z_i=\bar X_i/\bar Y_i$ & $C_i^+$ & Signal \\
\midrule
1  & 250 & 23.398 & 23.929 & 0.978 & 0.000 & 0 \\
2  & 250 & 26.696 & 26.970 & 0.990 & 0.000 & 0 \\
3  & 250 & 26.731 & 25.869 & 1.033 & 0.019 & 0 \\
4  & 250 & 23.047 & 23.082 & 0.998 & 0.003 & 0 \\
5  & 250 & 26.092 & 25.366 & 1.029 & 0.018 & 0 \\
6  & 250 & 24.346 & 22.623 & 1.076 & 0.080 & 0 \\
7  & 500 & 54.246 & 53.882 & 1.007 & 0.073 & 0 \\
8  & 500 & 49.963 & 49.233 & 1.015 & 0.073 & 0 \\
9  & 500 & 50.224 & 50.408 & 0.996 & 0.055 & 0 \\
10 & 500 & 50.675 & 48.288 & 1.049 & 0.091 & 0 \\
11 & 500 & 52.035 & 47.730 & 1.090 & 0.167 & 0 \\
12 & 500 & 59.143 & 56.362 & 1.049 & 0.202 & 0 \\
\textbf{13} & \textbf{500} & \textbf{57.159} & \textbf{49.633} & \textbf{1.152} & \textbf{0.339} & \textbf{1} \\
\textbf{14} & \textbf{500} & \textbf{58.934} & \textbf{50.854} & \textbf{1.159} & \textbf{0.484} & \textbf{1} \\
\textbf{15} & \textbf{250} & \textbf{28.036} & \textbf{25.631} & \textbf{1.094} & \textbf{0.564} & \textbf{1} \\
\bottomrule
\end{tabular}}
\end{table}

\section{Conclusions}
\label{sec:conclusions}

This paper proposed an upper one-sided adaptive cumulative sum chart for monitoring increases in the ratio of two correlated normal variables in short production runs. Unlike existing finite-horizon CUSUM-RZ schemes, in which the reference value $k$ is fixed in advance from a target shift size and only the decision interval $h$ is calibrated, the proposed chart jointly determines $(k^\star,h^\star)$ through a constrained bilevel optimization. The inner step calibrates $h(k)$ so that the in-control truncated average run length matches a prescribed target, and the outer step minimizes the out-of-control TARL$_1$ at a user-specified shift $\tau_{\mathrm{target}}$ over an admissible range of reference values. Both steps rely on a finite-state Markov-chain framework and an accurate approximation to the distribution of the subgroup ratio statistic.

The methodological advantages of the adaptive design are twofold. First, it yields a uniformly smaller (or equal) out-of-control TARL$_1$ than the fixed-$k$ design at the prescribed $\tau_{\mathrm{target}}$, with the largest gains arising in unequal coefficient-of-variation settings. Second, by allowing the inner step to search over reference values close to the in-control ratio $z_0$, the adaptive chart automatically resolves the boundary configurations of the fixed-$k$ design — those in which the prescribed TARL$_0$ cannot be attained because the upward signal probability under in-control becomes vanishingly small. This eliminates a documented numerical limitation of the recent short-run CUSUM-RZ literature.

The numerical analysis showed that the adaptive chart is highly effective for detecting moderate and large upward shifts in the ratio. In the matched-horizon benchmark, the adaptive CUSUM-RZ$^+$ chart performed on par with the EWMA-RZ and fixed-$k$ CUSUM-RZ$^+$ charts under a stable correlation structure and improved on both at the target shift when the correlation increased from Phase~I to Phase~II; all three memory-type charts substantially outperformed the short-run Shewhart-RZ chart of \citet{Tran2021}. The robustness analysis indicated that the chart is essentially insensitive to symmetric heavy-tailed departures from joint normality (Student-$t$ marginals with as few as $df=5$) but exhibits mild anti-conservative behaviour under contamination, with the achieved TARL$_0$ falling by approximately $5\%$--$11\%$ at $5\%$--$10\%$ contamination levels. The Phase~I study (200 replications per configuration) showed that a sample of $m\geq 100$ subgroups keeps the relative bias of the achieved TARL$_0$ at about $1\%$ of the target, while $m=25$ is associated with substantial run-length variability.

Several directions for future research remain. First, lower one-sided and two-sided versions of the adaptive chart could be developed for situations where the out-of-control direction is unknown or both directions are of practical concern. Second, the sensitivity to contamination motivates the use of robust estimators in the charting statistic, possibly along the lines of \citet{LucasCrosier1982} or \citet{HanQiao2023}. Third, a self-starting variant that updates $(k^\star,h^\star)$ with the cumulative information from Phase~II observations would alleviate the dependence on Phase~I sample size. Fourth, comparisons against the variable-sampling-interval CUSUM-RZ schemes of \citet{Yang2025} would be useful, although those require an additional layer of sampling-policy calibration that is beyond the scope of the present paper. Finally, the present framework can be extended to multivariate ratio structures and to settings with measurement error or autocorrelation, which are common in continuous-process industries.

\section*{Data availability statement}

No external empirical dataset was used in this study. The results are based on analytical calculations and numerical simulations. The simulated data and Python code used to reproduce the findings are available from the corresponding author upon reasonable request.

\section*{Declaration of Competing Interest}

The author declares no known competing financial interests or personal relationships that could have influenced the work reported in this paper.

\bibliographystyle{apacite}
\bibliography{paper}

\end{document}